\documentclass{amsart}
\usepackage{hyperref, url}
 \newtheorem{theorem}{Theorem}[section]
 
 \newtheorem{lemma}[theorem]{Lemma}

 \theoremstyle{definition}
 
 \theoremstyle{remark}

 \numberwithin{equation}{section}

\newcommand{\C}{\mathbb{C}}
\newcommand{\R}{\mathbb{R}}

\newcommand{\Z}{\mathbb{Z}}
\newcommand{\cA}{\mathcal{A}}
\newcommand{\cB}{\mathcal{B}}
\newcommand{\cC}{\mathcal{C}}
\newcommand{\cG}{\mathcal{G}}
\newcommand{\cJ}{\mathcal{J}}
\newcommand{\cK}{\mathcal{K}}
\newcommand{\cL}{\mathcal{L}}
\newcommand{\cM}{\mathcal{M}}
\newcommand{\cS}{\mathcal{S}}
\newcommand{\cZ}{\mathcal{Z}}
\newcommand{\alg}{\operatorname{alg}}
\newcommand{\diag}{\operatorname{diag}}
\newcommand{\im}{\operatorname{Im}}
\newcommand{\Ker}{\operatorname{Ker}}
\newcommand{\eps}{\varepsilon}
\begin{document}
%
%
%
%
%
%
%
%
%
\title[Singular Integral Operators]
{Singular Integral Operators\\
on Variable Lebesgue Spaces\\
over Arbitrary Carleson Curves}
\author[Karlovich]{Alexei Yu. Karlovich}
\address{
Departamento de Matem\'atica\\
Faculdade de Ci\^encias e Tecnologia\\
Universidade Nova de Lisboa\\
Quinta da Torre\\
2829--516 Caparica\\
Portugal}
\email{oyk@fct.unl.pt}
\thanks{The author is partially supported by the grant FCT/FEDER/POCTI/MAT/59972/2004.}


\dedicatory{To Professor Israel Gohberg on the occasion of his eightieth birthday}
\begin{abstract}
In 1968, Israel Gohberg and Naum Krupnik discovered that local spectra
of singular integral operators with piecewise continuous coefficients on
Lebesgue spaces $L^p(\Gamma)$ over Lyapunov curves have the shape of circular
arcs. About 25 years later, Albrecht B\"ottcher and Yuri Karlovich realized
that these circular arcs metamorphose to so-called logarithmic leaves with
a median separating point when Lyapunov curves metamorphose to arbitrary
Carleson curves. We show that this result remains valid in a more general
setting of variable Lebesgue spaces $L^{p(\cdot)}(\Gamma)$ where
$p:\Gamma\to(1,\infty)$ satisfies the Dini-Lipschitz condition. One of the
main ingredients of the proof is a new sufficient condition for the boundedness
of the Cauchy singular integral operator on variable Lebesgue spaces with
weights related to oscillations of Carleson curves.
\end{abstract}
\maketitle
\section{Introduction}
About forty years ago I.~Gohberg and N.~Krupnik \cite{GK68} constructed an
elegant Fredholm theory for singular integral operators with piecewise
continuous coefficients on Lebesgue spaces $L^p(\Gamma)$ over Lyapunov curves.
Their result says that the local spectra at discontinuity points of
the coefficients have the shape of circular arcs depending on $p$. That paper
was the starting point for generalizations and extensions of those results
to the case of power weights, matrix coefficients, and Banach algebras generated
by singular integral operators (see \cite{GK71,GK92}). I. Spitkovsky
\cite{Spitkovsky92} discovered that circular arcs metamorphose to massive horns
if one replaces power weights by general Muckenhoupt weights. A. B\"ottcher and
Yu. Karlovich \cite{BK95} observed that local spectra of singular integral
operators with piecewise continuous coefficients can be massive even on
$L^2(\Gamma)$ when $\Gamma$ is an arbitrary Carleson curve. The Fredholm theory
for the Banach algebra generated by matrix singular integral operators on
$L^p(\Gamma,w)$ under the most general conditions on the curve $\Gamma$ and
the weight $w$ is constructed by A.~B\"ottcher and Yu.~Karlovich and is
presented in the monograph \cite{BK97} (although, we advise to start the
study of this theory from the nice survey \cite{BK01}).

I.~Gohberg and N.~Krupnik \cite{GK68} also obtained some sufficient conditions
for the Fredholmness of singular integral operators with piecewise continuous
coefficients on so-called symmetric spaces (see \cite{KPS82} for the definition)
known also as rearrangement-invariant spaces (see \cite{BS88}). These spaces
include classical Lebesgue, Orlicz, and Lorentz spaces. The author
\cite{Karlovich98,Karlovich02} proved a criterion for the Fredholmness
of singular integral operators on rearrangement-invariant spaces and
observed that a ``complicated" space may also cause massiveness of local spectra.

Another natural generalization of the standard Lebesgue space $L^p(\Gamma)$ is
a so-called variable Lebesgue space $L^{p(\cdot)}$ defined in terms of the integral
\[
\int_\Gamma |f(\tau)|^{p(\tau)}\,|d\tau|
\]
(see the next section for the definition). Here the exponent $p$ is a continuous
function on $\Gamma$. Notice that variable Lebesgue spaces are not
rearrangement-invariant. V. Kokilashvili and S. Samko \cite{KS03} extended
the results of \cite{GK68} to the setting of variable Lebesgue spaces
over Lyapunov curves. In this setting, the circular arc depends
on the value of the exponent $p(t)$ at a discontinuity point $t\in\Gamma$.
Later on, the author gradually extended results known for singular integral
operators with piecewise continuous coefficients on weighted standard Lebesgue
spaces (see \cite{BK97,GK92}) to the case of weighted variable Lebesgue spaces
(see \cite{Karlovich03,Karlovich06} for power weights and Lyapunov curves;
\cite{Karlovich05} for power weights and so-called logarithmic Carleson curves;
\cite{Karlovich07} for radial oscillating weights and logarithmic Carleson curves).

In this paper we construct a Fredholm theory for the Banach algebra of singular
integral operators with matrix piecewise continuous coefficients on (non-weighted)
variable Lebesgue space over arbitrary Carleson curves.  We suppose that the
variable exponent is little bit better than continuous and, roughly speaking,
show that local spectra at the points $t$ of discontinuities of coefficients
are so-called logarithmic leaves (with a median separating point) \cite[Section~7.5]{BK97}
depending on the spirality indices $\delta_t^-,\delta_t^+$ of the curve at $t$
and the value $p(t)$. So that we replace the constant exponent $p$ in the results
for $L^p(\Gamma)$ \cite{BK95} by the value $p(t)$ at each point. Let us explain
why this is not so easy as it sounds. The only known method for studying singular
integral operators with piecewise continuous coefficients over arbitrary Carleson
curves is based on the Wiener-Hopf factorization technique, which in turn requires
an information on the boundedness of the Cauchy singular integral operator on spaces
with special weights related to oscillations of Carleson curves. For logarithmic
Carleson curves this boundedness problem is reduced to the case of power weights
treated in \cite{KPS06}. However, for arbitrary Carleson curves this is not the
case, a more general boundedness result was needed (formulated as a conjecture
in \cite{Karlovich08}). This need is satisfied in the present paper by a combination of
two very recent results by V.~Kokilashvili and S.~Samko \cite{KS08} and the author
\cite{Karlovich08}.

Let us also note that for standard Lebesgue spaces over slowly oscillating
Carleson curves (in particular, logarithmic Carleson curves) there
exists another method for studying singular integral operators based
on the technique of Mellin pseudodifferential operators and limit operators
(see e.g. \cite{Rabinovich96,BKR96,BKR00} and the references therein). It allows one to study
not only piecewise continuous coefficients but also coefficients
admitting discontinuities of slowly oscillating type. In this connection note
that very recently V. Rabinovich and S. Samko \cite{RS08} have started to study
pseudodifferential operators in the setting of variable Lebesgue spaces.
However, it seems that the method based on the Mellin technique does not
allow to consider the case of arbitrary Carleson curves.

The paper is organized as follows. In Section~\ref{sec:results} we give necessary
definitions and formulate the main results:
1) the above mentioned sufficient condition
for the boundedness of the Cauchy singular integral operator on a variable Lebesgue
space with a weight related to oscillations of an arbitrary Carleson curve;
2) a Fredholm criterion for an individual singular integral operator with piecewise
continuous coefficients in the spirit of \cite{GK68} and \cite{BK95};
3) a symbol calculus for the Banach algebra of singular integral operators with
matrix piecewise continuous coefficients. Sections~\ref{sec:boundedness}--\ref{sec:symbol}
contain the proofs of the results 1)--3), respectively.
\section{Preliminaries and main results}\label{sec:results}
\subsection{Carleson curves}
Let $\Gamma$ be a Jordan curve, that is, a curve that is homeomorphic to a circle.
We suppose that $\Gamma$ is rectifiable. We equip $\Gamma$ with Lebesgue length
measure $|d\tau|$ and the counter-clockwise orientation. The \textit{Cauchy
singular integral} of $f\in L^1(\Gamma)$ is defined by
\[
(Sf)(t):=\lim_{R\to 0}\frac{1}{\pi i}\int_{\Gamma\setminus\Gamma(t,R)}
\frac{f(\tau)}{\tau-t}d\tau
\quad (t\in\Gamma),
\]
where $\Gamma(t,R):=\{\tau\in\Gamma:|\tau-t|<R\}$ for $R>0$.
David \cite{David84} (see also \cite[Theorem~4.17]{BK97}) proved that the
Cauchy singular integral generates the bounded operator $S$ on the Lebesgue
space $L^p(\Gamma)$, $1<p<\infty$, if and only if $\Gamma$ is a
\textit{Carleson} (\textit{Ahlfors-David regular}) \textit{curve}, that is,
\[
\sup_{t\in\Gamma}\sup_{R>0}\frac{|\Gamma(t,R)|}{R}<\infty,
\]
where $|\Omega|$ denotes the measure of a measurable set $\Omega\subset\Gamma$.
\subsection{Variable Lebesgue spaces with weights}
A measurable function $w:\Gamma\to[0,\infty]$ is referred to as a \textit{weight
function} or simply a \textit{weight} if $0<w(\tau)<\infty$ for almost all
$\tau\in\Gamma$. Suppose $p:\Gamma\to(1,\infty)$ is a continuous function.
Denote by $L^{p(\cdot)}(\Gamma,w)$ the set of all measurable complex-valued
functions $f$ on $\Gamma$ such that
\[
\int_\Gamma |f(\tau)w(\tau)/\lambda|^{p(\tau)}|d\tau|<\infty
\]
for some $\lambda=\lambda(f)>0$. This set becomes a Banach space when equipped
with the Luxemburg-Nakano norm
\[
\|f\|_{p(\cdot),w}:=\inf\left\{\lambda>0:
\int_\Gamma |f(\tau)w(\tau)/\lambda|^{p(\tau)}|d\tau|\le 1\right\}.
\]
If $p$ is constant, then $L^{p(\cdot)}(\Gamma,w)$ is nothing else than
the weighted Lebesgue space. Therefore, it is natural to refer to
$L^{p(\cdot)}(\Gamma,w)$ as a \textit{weighted generalized Lebesgue space
with variable exponent} or simply as a \textit{weighted variable Lebesgue
space}. This is a special case of Musielak-Orlicz spaces \cite{Musielak83}
(see also \cite{KR91}). Nakano \cite{Nakano50} considered these spaces
(without weights) as examples of so-called modular spaces, and sometimes
the spaces $L^{p(\cdot)}(\Gamma,w)$ are referred to as weighted Nakano
spaces. In the case $w\equiv 1$ we will simply write $L^{p(\cdot)}(\Gamma)$.
\subsection{Boundedness of the Cauchy singular integral operator}
Let us define the weight we are interested in.
Fix $t\in\Gamma$ and consider the function $\eta_{t}:\Gamma\setminus\{t\}\to(0,\infty)$
defined by
\[
\eta_{t}(\tau):=e^{-\arg(\tau-t)},
\]
where $\arg(\tau-t)$ denotes any continuous branch of the argument on $\Gamma\setminus\{t\}$.
For every $\gamma\in\C$, put
\[
\varphi_{t,\gamma}(\tau):=|(\tau-t)^\gamma|=
|\tau-t|^{{\rm Re}\,\gamma}\eta_t(\tau)^{{\rm Im}\,\gamma}
\quad (\tau\in\Gamma\setminus\{t\}).
\]
A. ~B\"ottcher and Yu.~Karlovich \cite{BK95} (see also \cite[Chap.~1]{BK97})
proved that if $\Gamma$ is a Carleson Jordan curve, then at each point $t\in\Gamma$,
the following limits exist:
\[
\delta_t^-:=\lim_{x\to 0}\frac{\log(W_t^0\eta_t)(x)}{\log x},
\quad
\delta_t^+:=\lim_{x\to\infty}\frac{\log(W_t^0\eta_t)(x)}{\log x},
\]
where
\[
(W_t^0\eta_t)(x)=\limsup_{R\to 0}
\left(
\max_{\{\tau\in\Gamma:|\tau-t|=xR\}}\eta_t(\tau)\big/
\min_{\{\tau\in\Gamma:|\tau-t|=R\}}\eta_t(\tau)
\right).
\]
Moreover,
\[
-\infty<\delta_t^-\le\delta_t^+<+\infty.
\]
These numbers are called the lower and upper spirality indices of the curve
$\Gamma$ at $t$. For piecewise smooth curves $\delta_t^-\equiv\delta_t^+\equiv 0$,
for curves behaving like a logarithmic spiral in a neighborhood of $t$, one
has $\delta_t^-=\delta_t^+\ne 0$. However, the class of Carleson curves is
much larger: for all real numbers $-\infty<\alpha<\beta<+\infty$ there is
a Carleson curve $\Gamma$ such that $\delta_t^-=\alpha$ and $\delta_t^+=\beta$
at some $t\in\Gamma$ (see \cite[Proposition~1.21]{BK97}). Put
\[
\alpha_t^0(x):=\min\{\delta_t^-x,\delta_t^+x\},
\quad
\beta_t^0(x):=\max\{\delta_t^-x,\delta_t^+x\}
\quad(x\in\R).
\]

We will always suppose that $p:\Gamma\to(1,\infty)$ is a continuous function
satisfying the Dini-Lipschitz condition on $\Gamma$, that is, there exists a
constant $C_p>0$ such that
\[
|p(\tau)-p(t)|\le\frac{C_p}{-\log|\tau-t|}
\]
for all $\tau,t\in\Gamma$ such that $|\tau-t|\le1/2$.

Our first main result is the following theorem.
\begin{theorem}\label{th:boundedness-S}
Let $\Gamma$ be a Carleson Jordan curve and $p:\Gamma\to(1,\infty)$ be a
continuous function satisfying the Dini-Lipschitz condition. If $t\in\Gamma$,
$\gamma\in\C$, and
\begin{equation}\label{eq:boundedness-conditions}
0<\frac{1}{p(t)}+{\rm Re}\,\gamma+\alpha_t^0({\rm Im}\,\gamma),
\quad
\frac{1}{p(t)}+{\rm Re}\,\gamma+\beta_t^0({\rm Im}\,\gamma)<1,
\end{equation}
then $S$ is bounded on $L^{p(\cdot)}(\Gamma,\varphi_{t,\gamma})$.
\end{theorem}
For constant $p\in(1,\infty)$ this result (and the converse) is actually
proved by A.~B\"ottcher and Yu.~Karlovich \cite{BK95}, see also \cite{BK97}.
For a variable Lebesgue space with a power weight, that is, in the case when
${\rm Im}\,\gamma=0$, this result (and the converse) is due to V.~Kokilashvili,
V.~Paatashvili, and S.~Samko \cite{KPS06}. Note that V.~Kokilashvili, N.~Samko,
and S.~Samko \cite{KSS07-S} generalized that result also to the case of
so-called radial oscillating weights  $w_t(\tau)=f(|\tau-t|)$, where $f$ is
an oscillating function at zero. Obviously, $\eta_t$ is not of this type, in
general.

The proof of Theorem~\ref{th:boundedness-S} will be given in
Section~\ref{sec:boundedness}.
\subsection{Fredholm criterion}
Let $X$ be a Banach space and $\cB(X)$ be the Banach algebra of all bounded
linear operators on $X$. An operator $A\in\cB(X)$ is said to be Fredholm if
its image $\im A$ is closed in $X$ and the defect numbers
$\dim\Ker A$,  $\dim\Ker A^*$ are finite.
By $PC(\Gamma)$ we denote the set of all $a\in L^\infty(\Gamma)$ for which
the one-sided limits
\[
a(t\pm 0):=\lim_{\tau\to t\pm 0}a(\tau)
\]
exist and finite at each point $t\in\Gamma$; here $\tau\to t-0$ means that $\tau$
approaches $t$ following the orientation of $\Gamma$, while $\tau\to t+0$
means that $\tau$ goes to $t$ in the opposite direction. Functions in $PC(\Gamma)$
are called \textit{piecewise continuous} functions. Put
\[
P:=(I+S)/2,\quad Q:=(I-S)/2.
\]

By using Theorem~\ref{th:boundedness-S} and the machinery developed in
\cite{Karlovich03} (see also \cite{BK97}), we will prove our second main result.
\begin{theorem}\label{th:Fredholmness}
Let $\Gamma$ be a Carleson Jordan curve and $p:\Gamma\to(1,\infty)$ be a
continuous function satisfying the Dini-Lipschitz condition. Suppose
$a\in PC(\Gamma)$. The operator $aP+Q$ is Fredholm on $L^{p(\cdot)}(\Gamma)$
if and only if $a(t\pm 0)\ne 0$ and
\[
\frac{1}{p(t)}
-
\frac{1}{2\pi}\arg\frac{a(t-0)}{a(t+0)}+
\theta\alpha_t^0\left(\frac{1}{2\pi}\log\left|\frac{a(t-0)}{a(t+0)}\right|\right)+
(1-\theta)
\beta_t^0\left(\frac{1}{2\pi}\log\left|\frac{a(t-0)}{a(t+0)}\right|\right)
\]
is not an integer number for all $t\in\Gamma$ and all $\theta\in[0,1]$.
\end{theorem}
It is well known that  $\alpha_t^0(x)\equiv\beta_t^0(x)\equiv 0$ if $\Gamma$
is piecewise smooth. For Lyapunov curves and constant $p$, Theorem~\ref{th:Fredholmness}
was obtained by I.~Gohberg and N.~Krupnik \cite{GK68} (see also \cite[Chap.~9]{GK92}),
it was extended to variable Lebesgue spaces over Lyapunov curves or Radon
curves without cusps by V. Kokilashvili and S. Samko \cite{KS03}. For arbitrary
Carleson curves and constant $p$, Theorem~\ref{th:Fredholmness} is due to
A.~B\"ottcher and Yu.~Karlovich \cite{BK95} (see also \cite[Chap.~7]{BK97}).

The proof of Theorem~\ref{sec:Fredholmness} is presented in
Section~\ref{sec:Fredholmness}. It is developed following the well-known scheme
(see \cite{Spitkovsky92}, \cite[Chap.~7]{BK97}, and also
\cite{Karlovich03,Karlovich05,Karlovich07}).
\subsection{Leaves with a median separating point}
Let $p\in(0,1)$ and $\alpha,\beta:\R\to\R$ be continuous functions such that
$\alpha$ is concave, $\beta$ is convex, $\alpha(x)\le\beta(x)$ for all
$x\in\R$, and $0<1/p+\alpha(0)\le 1/p+\beta(0)<1$. Put
\[
Y(p,\alpha,\beta):=
\left\{\gamma=x+iy\in\C:\frac{1}{p}+\alpha(x)\le y\le\frac{1}{p}+\beta(x)\right\},
\]
and given $z_1,z_2\in\C$, let
\[
\cL(z_1,z_2;p,\alpha,\beta):=
\big\{M_{z_1,z_2}(e^{2\pi\gamma}):\gamma\in Y(p,\alpha,\beta)\big\}\cup\{z_1,z_2\},
\]
where
\[
M_{z_1,z_2}(\zeta):=(z_2\zeta-z_1)/(\zeta-1)
\]
is the M\"obius transform.
The set $\cL(z_1,z_2;p,\alpha,\beta)$ is referred to as the \textit{leaf} about
(or between) $z_1$ and $z_2$ determined by $p,\alpha,\beta$.

If $\alpha(x)=\beta(x)=0$ for all $x\in\R$, then
$\cL(z_1,z_2;p,\alpha,\beta)$
is nothing else than the circular arc
\[
\cA(z_1,z_2;p):=
\left\{z\in\C\setminus\{z_1,z_2\}:\arg\frac{z-z_1}{z-z_2}\in\frac{2\pi}{p}+2\pi\Z\right\}
\cup\{z_1,z_2\}.
\]

H.~Widom \cite{Widom60} was the first who understood the importance of these arcs
in the spectral theory of singular integral operators (in the setting of $L^p(\R)$).
These arcs play a very important role in the Gohberg-Krupnik Fredholm theory for
singular integral operators with piecewise continuous coefficients over Lyapunov
curves (see \cite{GK71,GK92}).

Suppose that $\alpha(x)=\beta(x)=\delta x$ for all $x\in\R$, where $\delta\in\R$. Then
the leaf $\cL(z_1,z_2;p,\alpha,\beta)$ is nothing else than the logarithmic double spiral
\[
\begin{split}
\cS(z_1,z_2;p,\delta)
:=&
\left\{z\in\C\setminus\{z_1,z_2\}:\arg\frac{z-z_1}{z-z_2}-\delta\log\left|\frac{z-z_1}{z-z_2}\right|
\in \frac{2\pi}{p}+2\pi\Z\right\}
\\
&\cup\{z_1,z_2\}.
\end{split}
\]

These logarithmic double spirals appear in the Fredholm theory for singular
integral operators over logarithmic Carleson curves, that is, when the spirality
indices $\delta_t^-$ and $\delta_t^+$ coincide at every point $t\in\Gamma$
(see \cite{BK97} and also \cite{Karlovich05}).

Now let $\delta^-$, $\delta^+$ be real numbers such that $\delta^-\le\delta^+$.
Suppose that
\begin{equation}\label{eq:indicators}
\alpha(x)=\min\{\delta^- x,\delta^+x\},\quad
\beta(x)=\max\{\delta^-x,\delta^+ x\}\quad
(x\in\R).
\end{equation}
Then it is not difficult to show that
\[
\cL(z_1,z_2;p,\alpha,\beta)=\bigcup_{\delta\in[\delta^-,\delta^+]}\cS(z_1,z_2;p,\delta).
\]
This set is always bounded by pieces of at most four logarithmic double spirals.
The point $m:=M_{z_1,z_2}(e^{2\pi i/p})$ has two interesting properties:
$m$ disconnects (separates) the leaf, that is, $\cL(z_1,z_2;p,\alpha,\beta)$
is connected, while $\cL(z_1,z_2;p,\alpha,\beta)\setminus\{m\}$ is a disconnected set;
and $m$ is a median point, that is, $|m-z_1|=|m-z_2|$ (see \cite[Example~7.10]{BK97}).
Leaves generated by functions $\alpha$ and $\beta$ of the form \eqref{eq:indicators}
are called logarithmic leaves with a median separating point. We refer to
\cite[Chap.~7]{BK97} for many nice plots of leaves (not only generated by
\eqref{eq:indicators}, but also more general).
\subsection{Symbol calculus for the Banach algebra of singular integral operators}
Let $N$ be a positive integer. We denote by $L_N^{p(\cdot)}$ the direct sum of
$N$ copies of $L^{p(\cdot)}(\Gamma)$ with the norm
\[
\|f\|=\|(f_1,\dots,f_N)\|:=(\|f_1\|_{p(\cdot)}^2+\dots+\|f_N\|_{p(\cdot)}^2)^{1/2}.
\]
The operator $S$ is defined on $L_N^{p(\cdot)}(\Gamma)$ elementwise.
We let stand $PC_{N\times N}(\Gamma)$ for the algebra of all $N\times N$
matrix functions with entries in $PC(\Gamma)$. Writing the elements of
$L_N^{p(\cdot)}(\Gamma)$ as columns, we can define the multiplication
operator $aI$ for $a\in PC_{N\times N}(\Gamma)$ as multiplication by the
matrix function $a$. Let $\cB:=\cB(L_N^{p(\cdot)}(\Gamma))$ be the Banach
algebra of all bounded linear operators on $L_N^{p(\cdot)}(\Gamma)$
and $\cK:=\cK(L_N^{p(\cdot)}(\Gamma))$ be its two-sided ideal consisting
of all compact operators on $L_N^{p(\cdot)}(\Gamma)$. By $\cA$ denote
the smallest closed subalgebra of $\cB$ containing the operator $S$ and
the set $\{aI:a\in PC_{N\times N}(\Gamma)\}$.

Our last main result is the following.
\begin{theorem}\label{th:symbol}
Suppose $\Gamma$ is a Carleson Jordan curve and $p:\Gamma\to(1,\infty)$
is a continuous function satisfying the Dini-Lipschitz condition.
Define the ``bundle" of logarithmic leaves with a median separating point by
\[
\cM:=\bigcup_{t\in\Gamma}\big(\{t\}\times
\cL(0,1;p(t),\alpha_t^0,\beta_t^0)\big).
\]
\begin{enumerate}
\item[{\rm(a)}] We have $\cK\subset\cA$.
\item[{\rm(b)}]
For each point $(t,z)\in\cM$, the map
\[
\sigma_{t,z}:\{S\}\cup\{aI\ :\ a\in PC_{N\times N}(\Gamma)\}\to\C^{2N\times 2N}
\]
given by
\[
\sigma_{t,z}(S)
=
\left[\begin{array}{cc}
E & O \\ O & -E
\end{array}\right],
\quad
\sigma_{t,z}(aI)
=
\]
\[
=\left[\begin{array}{cc}
a(t+0)z+a(t-0)(1-z) & (a(t+0)-a(t-0))\sqrt{z(1-z)}
\\[2mm]
(a(t+0)-a(t-0))\sqrt{z(1-z)} & a(t+0)(1-z)+a(t-0)z
\end{array}\right],
\]
where $E$ and $O$ denote the $N \times N$ identity and zero matrices,
respectively, and $\sqrt{z(1-z)}$ denotes any complex number whose square
is $z(1-z)$, extends to a Banach algebra homomorphism
\[
\sigma_{t,z}:\cA\to\C^{2N\times 2N}
\]
with the property that $\sigma_{t,z}(K)$ is the $2N\times 2N$ zero matrix
whenever $K$ is a compact operator on $L_N^{p(\cdot)}(\Gamma)$.

\item[{\rm(c)}]
An operator $A\in\cA$ is Fredholm on $L_N^{p(\cdot)}(\Gamma)$ if and only if
\[
\det\sigma_{t,z}(A)\ne 0
\quad\mbox{for all}\quad (t,z)\in\cM.
\]

\item[{\rm(d)}]
The quotient algebra $\cA/\cK$ is inverse closed in the Calkin algebra $\cB/\cK$,
that is, if a coset $A+\cK\in\cA/\cK$ is invertible in $\cB/\cK$, then
$(A+\cK)^{-1}\in\cA/\cK$.
\end{enumerate}
\end{theorem}
This theorem was proved by I.~Gohberg and N.~Krupnik for constant $p$
and Lyapunov curves (and power weights) in \cite{GK71}, and extended to the
setting of variable Lebesgue spaces over Lyapunov curves (again with power
weights) by the author \cite{Karlovich06}. The case of constant $p$
and arbitrary Carleson curves was treated by A.~B\"ottcher and Yu.~Karlovich
\cite{BK95} and Theorem~\ref{th:symbol} is a direct generalization of their
result to the setting of variable Lebesgue spaces.

We will present a sketch of the proof of Theorem~\ref{th:symbol}
in Section~\ref{sec:symbol}.
\section{Proof of the boundedness result}\label{sec:boundedness}
\subsection{Main ingredients}
It is well known that the boundedness of the Cauchy singular integral operator $S$
is closely related to the boundedness of the following maximal function
\[
(Mf)(t):=\sup_{\eps>0}\frac{1}{|\Gamma(t,\eps)|}\int_{\Gamma(t,\eps)}|f(\tau)|\,|d\tau|
\quad
(t\in\Gamma)
\]
defined for (locally) integrable functions $f$ on $\Gamma$. In particular, both
operators are bounded on weighted standard Lebesgue spaces $L^p(\Gamma,w)$
$(1<p<\infty)$ simultaneously  and this happen if and only if $w$ is a Muckenhoupt
weight. For weighted variable Lebesgue spaces a characterization of this sort
is unknown.

One of the main ingredients of the proof of Theorem~\ref{th:boundedness-S}
is the following very recent result by V.~Kokilashvili and S.~Samko.
\begin{theorem}[{\cite[Theorem~4.21]{KS08}}]
\label{th:KS}
Let $\Gamma$ be a Carleson Jordan curve. Suppose that $p:\Gamma\to(1,\infty)$
is a continuous function satisfying the Dini-Lipschitz condition and
$w:\Gamma\to[1,\infty]$ is a weight. If there exists a number $p_0$ such that
\[
1<p_0<\min\limits_{\tau\in\Gamma}p(\tau)
\]
and $M$ is bounded on $L^{p(\cdot)/(p(\cdot)-p_0)}(\Gamma,w^{-p_0})$, then
$S$ is bounded on $L^{p(\cdot)}(\Gamma,w)$.
\end{theorem}
The above conditional result allows us to derive sufficient conditions for
the boundedness of the Cauchy singular integral operator on weighted variable
Lebesgue spaces when some conditions for the boundedness of the maximal
operator in weighted variable Lebesgue spaces are already known.

Very recently, sufficient conditions for the boundedness of the maximal
operator that fit our needs are obtained by the author \cite{Karlovich08}
following the approach of \cite{KSS07-M}.
\begin{theorem}[{\cite[Theorem~4]{Karlovich08}}]
\label{th:boundedness-M}
Let $\Gamma$ be a Carleson Jordan curve. Suppose that $p:\Gamma\to(1,\infty)$
is a continuous function satisfying the Dini-Lipschitz condition. If $t\in\Gamma$,
$\gamma\in\C$, and \eqref{eq:boundedness-conditions} is fulfilled, then $M$
is bounded on $L^{p(\cdot)}(\Gamma,\varphi_{t,\gamma})$.
\end{theorem}
It should be noted that the author was not aware of \cite{KS08} at the time
of writing \cite{Karlovich08}. Because of that Theorem~\ref{th:boundedness-S}
was stated as a conjecture in \cite{Karlovich08}. Below we give its proof.
\subsection{Proof of Theorem~\ref{th:boundedness-S}}
Since $p:\Gamma\to(1,\infty)$ is continuous and $\Gamma$ is compact, we deduce
that $\min\limits_{\tau\in\Gamma}p(\tau)>1$. If the inequality
\[
\frac{1}{p(t)}+{\rm Re}\,\gamma+
\beta_t^0({\rm Im}\,\gamma)<1
\]
is fulfilled, then there exists a number $p_0$ such that
\[
1<p_0<\min\limits_{\tau\in\Gamma}p(\tau),
\quad
\frac{1}{p(t)}+{\rm Re}\,\gamma+\beta_t^0({\rm Im}\,\gamma)<\frac{1}{p_0}.
\]
The latter inequality is equivalent to
\begin{equation}\label{eq:boundedness-1}
0<
1-\frac{p_0}{p(t)}-
p_0\big({\rm Re}\,\gamma+\beta_t^0({\rm Im}\,\gamma)\big)
=
\frac{p(t)-p_0}{p(t)}-p_0{\rm Re}\,\gamma+
\alpha_t^0(-p_0{\rm Im}\,\gamma).
\end{equation}
Analogously, the inequality
\[
0<
\frac{1}{p(t)}+{\rm Re}\,\gamma+\alpha_t^0({\rm Im}\,\gamma)
\]
is equivalent to
\begin{equation}\label{eq:boundedness-2}
1 >
1-\frac{p_0}{p(t)}-
p_0\big({\rm Re}\,\gamma+
\alpha_t^0({\rm Im}\,\gamma)\big)
=\frac{p(t)-p_0}{p(t)}-p_0{\rm Re}\,\gamma+
\beta_t^0(-p_0{\rm Im}\,\gamma).
\end{equation}
From the equality $\varphi_{t,-p_0\gamma}=\varphi_{t,\gamma}^{-p_0}$,
inequalities \eqref{eq:boundedness-1}--\eqref{eq:boundedness-2},
and Theorem~\ref{th:boundedness-M} it follows that the maximal operator $M$
is bounded on $L^{p(\cdot)/(p(\cdot)-p_0)}(\Gamma,\varphi_{t,\gamma}^{-p_0})$.
To finish the proof, it remains to apply Theorem~\ref{th:KS}.
\qed
\section{Proof of the Fredholm criterion for the operator $aP+Q$}\label{sec:Fredholmness}
\subsection{Local representatives}\label{sec:local-representative}
In this section we suppose that $\Gamma$ is a Carleson Jordan curve and
$p:\Gamma\to(1,\infty)$ is a continuous functions satisfying the Dini-Lipschitz
condition.  Under these assumptions, the operator $S$ is bounded on
$L^{p(\cdot)}(\Gamma)$ by Theorem~\ref{th:boundedness-S}.

Functions $a,b\in L^\infty(\Gamma)$ are said to be locally equivalent at
a point $t\in\Gamma$ if
\[
\inf\big\{\|(a-b)c\|_\infty\ :\ c\in C(\Gamma),\ c(t)=1\big\}=0.
\]
\begin{theorem}\label{th:local_principle}
Suppose  $a\in L^\infty(\Gamma)$ and for each $t\in\Gamma$ there exists a function
$a_t\in L^\infty(\Gamma)$ which is locally equivalent to $a$ at $t$. If the
operators $a_tP+Q$ are Fredholm on $L^{p(\cdot)}(\Gamma)$ for all $t\in\Gamma$,
then $aP+Q$ is Fredholm on $L^{p(\cdot)}(\Gamma)$.
\end{theorem}
For weighted Lebesgue spaces this theorem is known as Simonenko's local
principle \cite{Simonenko65}. It follows from \cite[Theorem~6.13]{Karlovich03}.

The curve $\Gamma$ divides the complex plane $\mathbb{C}$ into the bounded
simply connected domain $D^+$ and the unbounded domain $D^-$. Without loss
of generality we assume that $0\in D^+$.
Fix $t\in\Gamma$.
For a function $a\in PC(\Gamma)$ such that $a^{-1}\in L^\infty(\Gamma)$,
we construct a ``canonical'' function $g_{t,\gamma}$ which is locally equivalent
to $a$ at the point $t\in\Gamma$. The interior and the exterior of the unit circle
can be conformally mapped onto $D^+$ and $D^-$ of $\Gamma$, respectively,
so that the point $1$ is mapped to $t$, and the points $0\in D^+$ and
$\infty\in D^-$ remain fixed. Let $\Lambda_0$ and $\Lambda_\infty$
denote the images of $[0,1]$ and $[1,\infty)\cup\{\infty\}$ under this map.
The curve $\Lambda_0\cup\Lambda_\infty$ joins $0$ to $\infty$ and
meets $\Gamma$ at exactly one point, namely $t$. Let $\arg z$ be a
continuous branch of argument in $\mathbb{C}\setminus(\Lambda_0\cup\Lambda_\infty)$.
For $\gamma\in\mathbb{C}$, define the function $z^\gamma:=|z|^\gamma e^{i\gamma\arg z}$,
where $z\in\mathbb{C}\setminus(\Lambda_0\cup\Lambda_\infty)$. Clearly, $z^\gamma$
is an analytic function in $\mathbb{C}\setminus(\Lambda_0\cup\Lambda_\infty)$. The
restriction of $z^\gamma$ to $\Gamma\setminus\{t\}$ will be denoted by
$g_{t,\gamma}$. Obviously, $g_{t,\gamma}$ is continuous and nonzero on
$\Gamma\setminus\{t\}$. Since $a(t\pm 0)\ne 0$, we can define
$\gamma_t=\gamma\in\mathbb{C}$ by the formulas
\begin{equation}\label{eq:local-representative}
\operatorname{Re}\gamma_t:=\frac{1}{2\pi}\arg\frac{a(t-0)}{a(t+0)},
\quad
\operatorname{Im}\gamma_t:=-\frac{1}{2\pi}\log\left|\frac{a(t-0)}{a(t+0)}\right|,
\end{equation}
where we can take any value of $\arg(a(t-0)/a(t+0))$, which implies that
any two choices of $\operatorname{Re}\gamma_t$ differ by an integer only.
Clearly, there is a constant $c_t\in\mathbb{C}\setminus\{0\}$ such that
$a(t\pm 0)=c_tg_{t,\gamma_t}(t\pm 0)$, which means that $a$ is locally
equivalent to $c_tg_{t,\gamma_t}$ at the point $t\in\Gamma$.
\subsection{Wiener-Hopf factorization of local representatives}
We say that a function $a\in L^\infty(\Gamma)$ admits a \textit{Wiener-Hopf
factorization on} $L^{p(\cdot)}(\Gamma)$ if $a^{-1}\in L^\infty(\Gamma)$ and
$a$ can be written in the form
\begin{equation}\label{eq:WH}
a(t)=a_-(t)t^\kappa a_+(t)
\quad\mbox{a.e. on}\ \Gamma,
\end{equation}
where $\kappa\in\Z$,  the factors $a_\pm$ enjoy the following properties:
\[
a_-\in QL^{p(\cdot)}(\Gamma)\stackrel{\cdot}{+}\mathbb{C},
\
a_-^{-1}\in QL^{q(\cdot)}(\Gamma)\stackrel{\cdot}{+}\mathbb{C},
\
a_+\in PL^{q(\cdot)}(\Gamma),
\
a_+^{-1}\in PL^{p(\cdot)}(\Gamma),
\]
where $1/p(t)+1/q(t)=1$ for all $t\in\Gamma$, and the operator $S$ is bounded
on the space $L^{p(\cdot)}(\Gamma,|a_+^{-1}|)$. One can prove that the number
$\kappa$ is uniquely determined.
\begin{theorem}\label{th:factorization}
A function $a\in L^\infty(\Gamma)$ admits a Wiener-Hopf factorization
\eqref{eq:WH} on $L^{p(\cdot)}(\Gamma)$ if and only if the operator $aP+Q$
is Fredholm on $L^{p(\cdot)}(\Gamma)$.
\end{theorem}
This theorem goes back to Simonenko \cite{Simonenko64,Simonenko68} for constant $p$.
For more about this topic we refer to \cite[Section~6.12]{BK97},
\cite[Section~5.5]{BS06}, \cite[Section~8.3]{GK92} and also to \cite{CG81,LS87}
in the case of weighted Lebesgue spaces. Theorem~\ref{th:factorization} follows
from \cite[Theorem~6.14]{Karlovich03}.

From \cite[Lemma~7.1]{Karlovich03} and the theorem on the boundedness of
the Cauchy singular integral operator on arbitrary Carleson curves
(see \cite{KPS06} or Theorem~\ref{th:boundedness-S}) we get the following
conditional result.
\begin{lemma}\label{le:fact-sufficiency}
If, for some $k\in\Z$ and $\gamma\in\C$, the operator $S$ is bounded on
the space $L^{p(\cdot)}(\Gamma,\varphi_{t,k-\gamma})$, then the function
$g_{t,\gamma}$ defined in Section~{\rm\ref{sec:local-representative}}
admits a Wiener-Hopf factorization on the space $L^{p(\cdot)}(\Gamma)$.
\end{lemma}
Combination of the above lemma and Theorem~\ref{th:boundedness-S} is the key
to the proof of the sufficiency portion of Theorem~\ref{th:Fredholmness}.
\subsection{Proof of Theorem~\ref{th:Fredholmness}}
\textit{Necessity.}
If $\Gamma$ is a Carleson Jordan curve, then $S$ is bounded on $L^{p(\cdot)}(\Gamma)$
(see \cite{KPS06} or Theorem~\ref{th:boundedness-S}). This implies that
the assumptions of \cite[Theorem~8.1]{Karlovich03} are satisfied.
Note that the indicator functions $\alpha_t$ and $\beta_t$ considered in
\cite[Theorem~8.1]{Karlovich03} (see also \cite[Chap.~3]{BK97}) coincide with
$\alpha_t^0$ and $\beta_t^0$, respectively, whenever we are in the nonweighted
situation (see e.g. \cite[Proposition~3.23]{BK97} or \cite[Lemma~5.15(a)]{Karlovich03}).
Therefore, the necessity portion of Theorem~\ref{th:Fredholmness} follows from
\cite[Theorem~8.1]{Karlovich03}.

\textit{Sufficiency.}
If $aP+Q$ is Fredholm on $L^{p(\cdot)}(\Gamma)$, then $a^{-1}\in L^\infty(\Gamma)$
in view of \cite[Theorem~6.11]{Karlovich03}. Therefore $a(t\pm 0)\ne 0$ for all
$t\in\Gamma$. Fix an arbitrary $t\in\Gamma$ and choose $\gamma=\gamma_t\in\C$
as in \eqref{eq:local-representative}. Then $a$ is locally equivalent to
$c_t g_{t,\gamma_t}$ at the point $t$, where $c_t$ is a nonzero constant
and the hypotheses of the theorem reads as follows:
\[
\frac{1}{p(t)}-{\rm Re}\,\gamma_t+
\theta\alpha_t^0(-{\rm Im}\,\gamma_t)+
(1-\theta)\beta_t^0(-{\rm Im}\,\gamma_t)\notin\Z
\mbox{ for all }\theta\in[0,1].
\]
Then there exists a number $k_t\in\Z$ such that
\[
0<\frac{1}{p(t)}+k_t-{\rm Re}\,\gamma_t+
\theta\alpha_t^0(-{\rm Im}\,\gamma_t)+
(1-\theta)\beta_t^0(-{\rm Im}\,\gamma_t)<1
\]
for all $\theta\in[0,1]$. In particular, if $\theta=1$, then
\begin{equation}\label{eq:sufficiency-1}
0<\frac{1}{p(t)}+{\rm Re}(k_t-\gamma_t)+\alpha_t^0({\rm Im}(k_t-\gamma_t));
\end{equation}
if $\theta=0$, then
\begin{equation}\label{eq:sufficiency-2}
\frac{1}{p(t)}+{\rm Re}(k_t-\gamma_t)+\beta_t^0({\rm Im}(k_t-\gamma_t))<1.
\end{equation}
From \eqref{eq:sufficiency-1}--\eqref{eq:sufficiency-2} and Theorem~\ref{th:boundedness-S}
it follows that the operator $S$ is bounded on $L^{p(\cdot)}(\Gamma,\varphi_{t,k_t-\gamma_t})$.
By Lemma~\ref{le:fact-sufficiency}, the function $g_{t,\gamma_t}$ admits a
Wiener-Hopf factorization on $L^{p(\cdot)}(\Gamma)$. Then, in view of
Theorem~\ref{th:factorization}, the operator $g_{t,\gamma_t}P+Q$ is Fredholm
on $L^{p(\cdot)}(\Gamma)$. It is easy to see that in this case the operator
$c_tg_{t,\gamma_t}P+Q$ is also Fredholm. Thus, for all local representatives
$c_tg_{t,\gamma_t}$, the operators $c_tg_{t\gamma_t}P+Q$ are Fredholm. To
finish the proof of the sufficiency part, it remains to apply the local principle
(Theorem~\ref{th:local_principle}).
\qed
\section{Construction of the symbol calculus}\label{sec:symbol}
\subsection{Allan-Douglas local principle}
In this section we present a sketch of the proof of Theorem~\ref{th:symbol}
based on the Allan-Douglas local principle and the two projections theorem
following the scheme of \cite[Chap.~8]{BK97} (see also
\cite{Karlovich05,Karlovich06,Karlovich07}).

Let $B$ be a Banach algebra with identity. A subalgebra $Z$ of $B$ is said to
be a \textit{central subalgebra} if $zb=bz$ for all $z\in Z$ and all $b\in B$.
\begin{theorem}[see \cite{BS06}, Theorem~1.35(a)]
\label{th:AllanDouglas}
Let $B$ be a Banach algebra with identity $e$ and let $Z$ be a closed central
subalgebra of $B$ containing $e$. Let $M(Z)$ be the maximal ideal space of $Z$,
and for $\omega\in M(Z)$, let $J_\omega$ refer to the smallest closed two-sided
ideal of $B$ containing the ideal $\omega$. Then an element $b$ is invertible
in $B$ if and only if $b+J_\omega$ is invertible in the quotient algebra
$B/J_\omega$ for all $\omega\in M(Z)$.
\end{theorem}
The algebra $B/J_\omega$ is referred to as the \textit{local algebra} of $B$
at $\omega\in M(Z)$ and the spectrum of $b+J_\omega$
in $B/J_\omega$ is called the \textit{local spectrum} of $b$ at $\omega\in M(Z)$.
\subsection{Localization}\label{sec:localization}
An operator $A\in\cB$ is said to be of local type if its commutator with
the operator of multiplication by the diagonal matrix function $\diag\{c,\dots,c\}$
is compact for every continuous function $c$ on $\Gamma$. The set $\cL$ of all
operators of local type forms a Banach subalgebra of $\cB$. By analogy with
\cite[Lemma~5.1]{Karlovich06} one can prove that $\cK\subset\cL$. From
\cite[Lemma~6.5]{Karlovich03} it follows that the operator $S$ is of local type.
Thus,
\[
\cK\subset\cA\subset\cL.
\]
It is easy to see that $A\in\cL$ is Fredholm if and only if the coset $A+\cK$
is invertible in $\cL/\cK$. We will study the invertibility of a coset
$A+\cK$ of $\cA/\cK$ in the larger algebra $\cL/\cK$ by using the Allan-Douglas
local principle. Consider
\[
\cZ/\cK:=\{\diag\{c,\dots,c\}I+\cK:c\in C(\Gamma)\}.
\]
Every element of this subalgebra commutes with all elements of $\cL/\cK$.
The maximal ideal spaces $M(\cZ/\cK)$ of $\cZ/\cK$ may be identified with
the curve $\Gamma$ via the Gelfand map
\[
\cG:\cZ/\cK\to C(\Gamma),
\quad
\big(\cG(\diag\{c,\dots,c\}I+\cK)\big)(t)=c(t)
\quad(t\in\Gamma).
\]
For every $t\in\Gamma$ we define $\cJ_t\subset\cL/\cK$ as the smallest
closed two-sided ideal of $\cL/\cK$ containing the set
\[
\big\{\diag\{c,\dots,c\}I+\cK:\ c\in C(\Gamma), \ c(t)=0\big\}.
\]
Let $\chi_t$ be the characteristic of a proper arc of $\Gamma$ starting at $t\in\Gamma$.
For a matrix function $a\in PC_{N\times N}(\Gamma)$, let
\[
a_t:=a(t-0)(1-\chi_t)+a(t+0)\chi_t.
\]
It is easy to see that $aI-a_tI+\cK\in\cJ_t$. This implies that for any operator
$A\in\cA$, the coset $A+\cK+\cJ_t$ belongs to the smallest closed subalgebra $\cA_t$
of the algebra $\cL_t:=(\cL/\cK)/\cJ_t$ that contains the cosets
\begin{equation}\label{eq:projections-p-and-q}
p:=P+\cK+\cJ_t,
\quad
q:=\diag\{\chi_t,\dots,\chi_t\}I+\cK+\cJ_t
\end{equation}
and the algebra
\begin{equation}\label{eq:algebra-C}
\cC:=\big\{cI+\cK+\cJ_t: \ c\in\C^{N\times N}\big\}.
\end{equation}
Thus, by the Allan-Douglas local principle,
for every $A\in\cA$, the problem of invertibility of $A+\cK$ in the
algebra $\cL/\cK$ is reduced to the problem of invertibility of
$A+\cK+\cJ_t\in\cA_t$ in the local algebra $\cL_t$ for every $t\in\Gamma$.
\subsection{The two projections theorem}
Recall that an element $r$ of a Banach algebra is called an \textit{idempotent}
(or, somewhat loosely, also a \textit{projection}), if $r^2=r$.

The following two projections theorem was obtained by T.~Finck,
S.~Roch, and B.~Silbermann \cite{FRS93} and in a slightly different form by
I.~Gohberg and  N.~Krupnik \cite{GK93}
(see also \cite[Section~8.3]{BK97}).
\begin{theorem}\label{th:2proj}
Let $B$ be a Banach algebra with identity $e$, let $\cC$ be a Banach subalgebra
of $B$ which contains $e$ and is isomorphic to ${\mathbb{C}}^{N \times N}$,
and let $r$ and $s$ be two idempotent elements in $B$ such that
$cr=rc$ and $cs=sc$ for all $c \in \cC$. Let $A=\alg(\cC,r,s)$
be the smallest closed subalgebra of $B$ containing $\cC,r,s$. Put
\[
x=rsr+(e-r)(e-s)(e-r),
\]
denote by $\mathrm{sp}\,x$ the spectrum of $x$ in $B$, and suppose the points
$0$ and $1$ are not isolated points of $\mathrm{sp}\,x$. Then
\begin{enumerate}
\item[{\rm (a)}]
for each $z \in \mathrm{sp}\,x$ the map $\sigma_{z}$ of $\cC \cup \{r,s\}$
into the algebra ${\mathbb{C}}^{2N\times 2N}$ of all complex $2N\times 2N$
matrices defined by
\begin{eqnarray*}
&&
\sigma_{z}c=\left[
\begin{array}{cc}
c & O\\
O & c
\end{array}
\right],
\sigma_{z}r=\left[
\begin{array}{cc}
E & O\\
O & O
\end{array}
\right],
\sigma_{z}s=\left[
\begin{array}{cc}
z E & \sqrt{z(1-z)}E \\
\sqrt{z(1-z)}E  & (1-z)E
\end{array}
\right],
\end{eqnarray*}
where $c\in \cC$, $E$ and $O$ denote the $N \times N$ identity and zero matrices,
respectively, and $\sqrt{z(1-z)}$ denotes any complex number whose square is $z(1-z)$,
extends to a Banach algebra homomorphism
\[
\sigma_{z}: A \to {\mathbb{C}}^{2N \times 2N};
\]

\item[{\rm (b)}]
every element $a$ of the algebra $A$ is invertible in the algebra $B$ if
and only if
\[
\det \sigma_{z} a \neq 0 \quad\mbox{for all}\quad z \in \mathrm{sp}\,x;
\]

\item[{\rm (c)}]
the algebra $A$ is inverse closed in $B$ if and only if the spectrum of
$x$ in $A$ coincides with the spectrum of $x$ in $B$.
\end{enumerate}
\end{theorem}
\subsection{Local algebras $\cA_t$ and $\cL_t$ are subject
to the two projections theorem}
In this subsection we verify that the algebras $\cA_t$ and $\cL_t$ defined
in Section~\ref{sec:localization} satisfy the assumptions of the two
projections theorem (Theorem~\ref{th:2proj}). It is obvious that the algebra
$\cC$ defined by \eqref{eq:algebra-C} is isomorphic to the algebra $\C^{N\times N}$.
It is easy to see also that
\[
p^2=p,\quad
q^2=q,\quad
pc=sp,\quad
qc=cq
\]
for all $c\in\cC$.

From Theorem~\ref{th:Fredholmness} by analogy with \cite[Theorem~8.19]{BK97}
one can derive the following.
\begin{theorem}
Let $t\in\Gamma$ and the elements $p,q\in\cA_t$ be given by
\eqref{eq:projections-p-and-q}. The spectrum of the element
\[
x:=pqp+(e-p)(e-q)(e-p)
\]
in the algebra $\cL_t$ coincides with the logarithmic leaf with
a median separating point $\cL(0,1;p(t),\alpha_t^0,\beta_t^0)$.
\end{theorem}

Notice that $0$ and $1$ are not isolated points of the leaf
$\cL(0,1;p(t),\alpha_t^0,\beta_t^0)$.

We have shown that $\cA_t$ and $\cL_t$ satisfy all the assumptions of
the two projections theorem. Thus, our last main result (Theorem~\ref{th:symbol})
is obtained by localizing as above and then by applying the two projections
theorem to the local algebras $\cA_t$ and $\cL_t$ (see \cite{BK97}
and also \cite{Karlovich05,Karlovich06,Karlovich07} for more details).
We only note that the mapping $\sigma_{t,z}$ in Theorem~\ref{th:symbol}
is constructed by the formula
\[
\sigma_{t,z}=\sigma_z\circ\pi_t,
\]
where $\sigma_z$ is the mapping from Theorem~\ref{th:2proj} and
$\pi_t$ acts by the rule $A\mapsto A+\cK+\cJ_t$.


\begin{thebibliography}{99}
\bibitem{BS88}
C. Bennett and R. Sharpley,
\textit{Interpolation of Operators.}
Academic Press, Boston, 1988.

\bibitem{BK95}
A. B\"ottcher and Yu. I. Karlovich,
\textit{Toeplitz and singular integral operators on Carleson curves with
logarithmic whirl points.}
\href{http://dx.doi.org/10.1007/BF01208347}
{Integral Equations Operator Theory \textbf{22} (1995), 127--161}.

\bibitem{BK97}
A. B\"ottcher and Yu. I. Karlovich,
\textit{Carleson Curves, Muckenhoupt Weights, and Toeplitz Operators.}
Birkh\"auser, Basel, 1997.

\bibitem{BK01}
A. B\"ottcher and Yu. I. Karlovich,
\textit{Cauchy's singular integral operator and its beautiful spectrum.}
In: ``Systems, approximation, singular integral operators,
and related topics'' (Bordeaux, 2000).
Operator Theory: Advances and Applications \textbf{129} (2001), 109--142.

\bibitem{BKR96}
A. B\"ottcher, Yu. I. Karlovich, and V. S. Rabinovich,
\textit{Emergence, persistence, and disappearance of logarithmic spirals in the
spectra of singular integral operators.}
\href{http://dx.doi.org/10.1007/BF01203026}
{Integral Equations Operator Theory \textbf{25} (1996), 406--444}.

\bibitem{BKR00}
A. B\"ottcher, Yu. I. Karlovich, and V. S. Rabinovich,
\textit{The method of limit operators for one-dimensional singular integrals
with slowly oscillating data.}
\href{http://www.mathjournals.org/jot/2000-043-001/2000-043-001-008.pdf}
{J. Operator Theory \textbf{43} (2000), 171--198}.

\bibitem{BS06}
A. B\"ottcher and B. Silbermann,
\href{http://dx.doi.org/10.1007/3-540-32436-4}
{\textit{Analysis of Toeplitz Operators}}.
2nd edition.
Springer-Verlag, Berlin, 2006.

\bibitem{CG81}
K. P. Clancey and I. Gohberg,
\textit{Factorization of Matrix Functions and Singular Integral Operators.}
Operator Theory: Advances and Applications \textbf{3}.
Birkh\"auser, Basel, 1981.

\bibitem{David84}
G. David,
\textit{Oper\'ateurs int\'egraux singuliers sur certaines courbes du plan
complexe.}
\href{http://www.numdam.org/item?id=ASENS_1984_4_17_1_157_0}
{Ann. Sci. \'Ecole Norm. Super. \textbf{17} (1984), 157--189}.

\bibitem{FRS93}
T. Finck, S. Roch, and B. Silbermann,
\textit{Two projections theorems and symbol calculus for operators
with massive local spectra.}
\href{http://dx.doi.org/10.1002/mana.19931620114}
{Math. Nachr. \textbf{162} (1993), 167--185}.

\bibitem{GK68}
I. C. Gohberg and N. Ya. Krupnik,
\textit{The spectrum of singular integral operators in $L\sb{p}$ spaces}.
\href{http://matwbn.icm.edu.pl/ksiazki/sm/sm31/sm31131.pdf}
{Studia Math. \textbf{31} (1968), 347--362}
(in Russian).

\bibitem{GK71}
I. Gohberg and N. Krupnik,
\textit{Singular integral operators with piecewise continuous coefficients
and their symbols.}
\href{http://dx.doi.org/10.1070/IM1971v005n04ABEH001127}
{Math. USSR Izvestiya \textbf{5} (1971), 955--979}.

\bibitem{GK92}
I. Gohberg and N. Krupnik,
\textit{One-Dimensional Linear Singular Integral Equations.}
Vols. 1 and 2.
Operator Theory: Advances and Applications \textbf{53--54}.
Birkh\"auser, Basel, 1992.

\bibitem{GK93}
I. Gohberg and N. Krupnik,
\textit{Extension theorems for Fredholm and invertibility symbols.}
\href{http://dx.doi.org/10.1007/BF01205291}
{Integral Equations Operator Theory \textbf{16} (1993), 514--529}.

\bibitem{Karlovich98}
A. Yu. Karlovich,
\textit{Singular integral operators with piecewise continuous coefficients
in reflexive rearrangement-invariant spaces.}
\href{http://dx.doi.org/10.1007/BF01194990}
{Integral Equatations Operator Theory \textbf{32} (1998), 436--481}.

\bibitem{Karlovich02}
A. Yu. Karlovich,
\textit{Algebras of singular integral operators with $PC$ coefficients
in rearrangement-invari\-ant spaces with Muckenhoupt weights.}
\href{http://www.mathjournals.org/jot/2002-047-002/2002-047-002-004.pdf}
{J. Operator Theory \textbf{47} (2002), 303--323}.

\bibitem{Karlovich03}
A. Yu. Karlovich,
\textit{Fredholmness of singular integral operators with piecewise continuous
coefficients on weighted Banach function spaces.}
\href{http://dx.doi.org/doi:10.1216/jiea/1181074970}
{J. Integr. Equat. Appl. \textbf{15} (2003), 263--320}.

\bibitem{Karlovich05}
A. Yu. Karlovich,
\textit{Algebras of singular integral operators on Nakano spaces with
Khve\-delidze weights over Carleson curves with logarithmic whirl points.}
In: ``Pseudodifferential Equations and Some Problems of
Mathematical Physics", Rostov-on-Don, 2005, 135--142.
Preprint is available at
\href{http://arxiv.org/abs/math/0507312}
{arXiv:math/0507312}.

\bibitem{Karlovich06}
A. Yu. Karlovich,
\textit{Algebras of singular integral operators with piecewise
continuous coefficients on weighted Nakano spaces.}
In: ``The Extended Field of Operator Theory".
\href{http://dx.doi.org/10.1007/978-3-7643-7980-3}
{Operator Theory: Advances and Applications \textbf{171} (2006), 171--188}.

\bibitem{Karlovich07}
A. Yu. Karlovich,
\textit{Singular integral operators on variable Lebesgue spaces with radial oscillating weights}.
Preprint, 2007,
\href{http://arxiv.org/abs/0708.0778}{arXiv:0708.0778}.

\bibitem{Karlovich08}
A. Yu. Karlovich,
\textit{Maximal operators on variable Lebesgue spaces with weights related to
oscillations of Carleson curves}.
Preprint, 2008,
\href{http://arxiv.org/abs/0808.0258}{arXiv:0808.0258}.

\bibitem{KPS06}
V. Kokilashvili, V. Paatashvili, and S. Samko,
\textit{Boundedness in Lebesgue spaces with variable exponent of the Cauchy
singular operator on Carleson curves.}
In: ``Modern Operator Theory and Applications. The Igor Borisovich Simonenko
Anniversary Volume".
\href{http://dx.doi.org/10.1007/978-3-7643-7737-3}
{Operator Theory: Advances and Applications \textbf{170} (2006), 167--186}.

\bibitem{KSS07-M}
V. Kokilashvili, N. Samko, and S. Samko,
\textit{The maximal operator in weighted variable spaces $L^{p(\cdot)}$}.
J. Funct. Spaces Appl. \textbf{5} (2007), 299--317.

\bibitem{KSS07-S}
V. Kokilashvili, N. Samko, and S. Samko,
\textit{Singular operators in variable spaces $L^{p(\cdot)}(\Omega,\rho)$
with oscillating weights.}
\href{http://dx.doi.org/10.1002/mana.200510542}
{Math. Nachr. \textbf{280} (2007), 1145--1156}.

\bibitem{KS03}
V. Kokilashvili and S. Samko,
\textit{Singular integral equations in the Lebesgue spaces with variable exponent.}
Proc. A. Razmadze Math. Inst. \textbf{131} (2003), 61--78.

\bibitem{KS08}
V. Kokilashvili and S. Samko,
\textit{Operators of harmonic analysis in weighted spaces with non-standard growth.}
\href{http://dx.doi.org/10.1016/j.jmaa.2008.06.056}
{J. Math. Anal. Appl., doi:10.1016/j.jmaa.2008.06.056}.

\bibitem{KR91}
O. Kov\'a{\v c}ik and J.~R\'akosn{\'\i}k,
\textit{On spaces $L\sp {p(x)}$ and $W\sp {k,p(x)}$.}
\href{http://hdl.handle.net/10338.dmlcz/102493}
{Czechoslovak Math. J. \textbf{41(116)} (1991), 592--618}.

\bibitem{KPS82}
S. G. Krein, Ju. I. Petunin, and E. M. Semenov,
\textit{Interpolation of Linear Operators.}
AMS Translations of Mathematical Monographs \textbf{54},
Providence, RI, 1982.

\bibitem{LS87}
G. S. Litvinchuk and I. M. Spitkovsky,
\textit{Factorization of Measurable Matrix Functions.}
Operator Theory: Advances and Applications \textbf{25}.
Birkh\"auser, Basel, 1987.

\bibitem{Musielak83}
J. Musielak,
\textit{Orlicz Spaces and Modular Spaces.}
\href{http://dx.doi.org/10.1007/BFb0072210}
{Lecture Notes in Mathematics \textbf{1034}}.
Springer-Verlag, Berlin, 1983.

\bibitem{Nakano50}
H. Nakano,
\textit{Modulared Semi-Ordered Linear Spaces.}
Maruzen Co., Ltd., Tokyo, 1950.

\bibitem{Rabinovich96}
V. S. Rabinovich,
\textit{Algebras of singular integral operators on composed contours with
nodes that are logarithmic whirl points.}
\href{http://dx.doi.org/10.1070/IM1996v060n06ABEH000099}
{Izvestiya Mathematics \textbf{60} (1996), 1261--1292}.

\bibitem{RS08}
V. Rabinovich and S. Samko,
\textit{Boundedness and Fredholmness of pseudodifferential operators in
variable exponent spaces.}
\href{http://dx.doi.org/10.1007/s00020-008-1566-9}
{Integral Equations Operator Theory \textbf{60} (2008), 507--537}.

\bibitem{Simonenko64}
I. B. Simonenko,
\textit{The Riemann boundary value problem for $n$ pairs functions with
measurable coefficients and its application to the investigation of
singular integral operators in the spaces $L^p$ with weight.}
Izv. Akad. Nauk SSSR, Ser. Matem. \textbf{28} (1964), 277--306 (in Russian).

\bibitem{Simonenko65}
I. B. Simonenko,
\textit{A new general method of investigating linear operator equations
of singular integral equations type.}
Part I:  Izv. Akad. Nauk SSSR, Ser. Matem. \textbf{29} (1965), 567--586
(in Russian);
Part II: Izv. Akad. Nauk SSSR, Ser. Matem. \textbf{29} (1965), 757--782
(in Russian).

\bibitem{Simonenko68}
I. B. Simonenko,
\textit{Some general questions in the theory of the Riemann boundary value problem.}
\href{http://dx.doi.org/10.1070/IM1968v002n05ABEH000706}
{Math. USSR Izvestiya \textbf{2} (1968), 1091--1099}.

\bibitem{Spitkovsky92}
I. M. Spitkovsky,
\textit{Singular integral operators with $PC$ symbols on the spaces with
general weights.}
\href{http://dx.doi.org/10.1016/0022-1236(92)90075-T}
{J. Funct. Anal. \textbf{105} (1992), 129--143}.

\bibitem{Widom60}
H. Widom,
\textit{Singular integral equations on $L^p$.}
\href{http://www.jstor.org/stable/1993367}
{Trans. Amer. Math. Soc. \textbf{97} (1960), 131--160}.
\end{thebibliography}
\end{document}